# ON RUGINA'S SYSTEM OF THOUGHT

by Florentin Smarandache, Ph. D.
Mathematics Department
University of New Mexico
200 College Road
Gallup, NM 87301.
E-mail: smarand@unm.edu

**Abstract**. In this article one investigates Rugina's Orientation Table and one gives particular examples for several of its seven models.
Leon Walras's Economics of Stable Equilibrium and Keynes's Economics of Disequilibrium are combined in Rugina's Orientation Table in systems which are s% stable and 100-s% unstable, where s may be 100, 95, 65, 50, 35, 5, and 0.
The Classical Logic and Modern Logic are united in Rugina's Integrated Logic, and then generalized in the Neutrosophic Logic.

## I) Introduction.

Coming across Rugina's System of Thoughts, in his published books and articles [7-11], I learned about the *connection* between classic and modern. It is not a contradiction, but a complementarity from the part of modern with respect to the classic; and always the new 'modern' will have something to bring to the old knowledge.
In a similar way we may talk on the complementarity between theory and practices, rather than their contradiction.
In economics, Rugina negated Marx's social justice for the mass and Keynes's involuntary unemployment. His methodology in science tries to unite all scientific fields, preserving however independence in thinking and judgement.

Einstein worked in the last period of his life on the Unified Field Theory (a single general theory in physics), but didn't succeed. At the present, his supposition that the speed of light is a barrier in the universe is also being denied.

The economical systems are characterized by free market or centrally-planned and controlled economy. I think each system has a mixture of the previous, where a part of the market is free and another is centrally-planned and controlled.

Rugina's Universal Hypothesis of Duality:
The physical universe is composed of stable and unstable elements arranged in various proportions,
may be completed with unknown elements, a strip border between stable and unstable, which are continously changing from the state of equilibrium to disequilibrium and vice-versa, and which therefore are giving the dynamics of the universe.
Unknown may be: anomalies, relativities, uncertainties, revolution risks, hidden parameters.
The internal parameters are involved in Rugina's Universal Law of Natural Parameter (NaPa):
Any system in order to reach and maintain a position of stable equilibrium must have a very strong natural parameter (center of weight).
Whereas the external parameters are involved in Rugina's Universal Law of General Consistency:
Any system produces and maintains a position of stable equilibrium if there is a suitable space-time frame work.

## II) Theory of Paradoxes.

How did I get to the Theory of Paradoxes?
I have observed that: what's good for some ones, may be bad for others - and reciprocally. There are peoples who are considered terrorists by their enemies, and patriots by their friends. All of them are right and wrong in the same time. If one changes the referential system, the result is different.

Let's see a few nice paradoxes:

**Social Paradox**:
In a democracy should the nondemocratic ideas be allowed?

a) If yes, i.e. the nondemocratic ideas are allowed, then one not has a democracy anymore. (The nondemocratic ideas may overturn the society.)

b) If no, i.e. other ideas are not allowed - even those nondemocratic -, then one not has a democracy either, because the freedom of speech is restricted.

**The Sets' Paradox**:
The notion of "set of all sets", introduced by Georg Cantor, does not exist.
Let all sets be noted by $\{S_a\}_a$, where a indexes them.
But the set of all sets is itself another set, say $T_1$;
and then one constructs again another "set of <all sets>", but <all sets> are this time $\{S_a\}$ and $T_1$, and then the "set of all sets" is now $T_2$, different from $T_1$;
and so on... .
Even the notion of "all sets" can not exactly be defined (like the largest number of an open interval, which doesn't exist), as one was just seeing above (we can construct a new set as the "set of all sets") and reunites it to "all sets".

A **Paradoxist Psychological Cómplex** (with the accent on the first syllable):

A collection of fears stemming from previous unsuccessful experience or from unconscious feelings that, wanting to do something <S>, the result would be <Anti-S>,
which give rise to feelings, attitudes, and ideas pushing the subject towards a deviation of action <S> eventually towards an <Anti-S> action.
(From the positive and negative brain's electrical activities.)
For example: A shy boy, attempting to invite a girl to dance, inhabits himself of fear she would turn him down...

How to manage this phobia? To dote and anti-dote!
By transforming it into an opposite one, thinking differently, and being fear in our mind that we would pass our expectancies but we shouldn't.

People who do not try of fear not to be rejected: they lose by not competing!

**Auto-suggestion**:
If an army leaves for war with anxiety to lose, that army are half-defeated before starting the confrontation.

**Paradoxist Psychological Behavior**:
How can we explain contrary behaviors of a person: in the same conditions, without any reason, cause?
Because our deep unconsciousness is formed of contraries.

**Ceaseless Anxiety**:
What you want is, normally, what you don't get. And this is for eternity. Like a chain...
Because, when you get it (if ever), something else will be your next desire. Man can't live without a new hope.

**Inverse Desire**:
The wish to purposely have bad luck, to suffer, to be pessimistic as stimulating factors for more and better creation or work.

(Applies to some artists, poets, painters, sculptors, spiritualists.)

All is possible, the impossible too!
Is this an optimistic or pessimistic paradox?
a) It is an optimistic paradox, because shows that all is possible.
b) It is a pessimistic paradox, because shows that the impossible is possible.

**Mathematician's Paradox**:
Let M be a mathematician who may not be characterized by his mathematical work.
a) To be a mathematician, M should have some mathematical work done, therefore M should be characterized by that work.
b) The reverse judgement: if M may not be characterized by his mathematical work, then M is not a mathematician.

**Divine Paradox (I)**:
Can God commit suicide?
If God cannot, then it appears that there is something God cannot do, therefore God is not omnipotent.
If God can commit suicide, then God dies - because He has to prove it, therefore God is not immortal.

**Divine Paradox (II)**:
Can God be atheist, governed by scientific laws?
If God can be atheist, then God doesn't believe in Himself, therefore why should we believe in Him?
If God cannot, then again He's not omnipotent.

[Religion is full of god-ism and evil-ism.]

**God and Evil in the same Being**.
Man is a bearer of good and bad simultaneously. Man is enemy to himself. God and Magog!

**Expect the Unexpected**:
If we expect someone to do the unexpected, then:
- is it possible for him to do the unexpected?
- is it possible for him to do the expected?
If he does the unexpected, then that's what we expected.
If he doesn't do the expected, then he did the unexpected.

**The Ultimate Paradox**:
Living is the process of dying.
Reciprocally: Death of one is the process of somebody else's life [an animal eating another one].

Exercises for readers:
If China and Japan are in the Far East, why from USA do we go west to get there?
Are humans inhuman, because they committed genocides?

**The Invisible Paradoxes**:
Our visible world is composed of a totality of invisible particles.
Things with mass result from atoms with quasi-null mass.
Infinity is formed of finite part(icle)s.
Look at these Sorites Paradoxes (associated with Eubulides of Miletus (fourth century B.C.):

a) An invisible particle does not form a visible object, nor do two invisible particles, three invisible particles, etc. However, at some point, the collection of invisible particles becomes large enough to form a visible object, but there is apparently no definite point where this occurs.

b) A similar paradox is developed in an opposite direction. It is always possible to remove an atom from an object in such a way that what is left is still a visible object. However, repeating and repeating this process, at some point, the visible object is decomposed so that the left part becomes invisible, but there is no definite point where this occurs.
Between <A> and <Non-A> there is no clear distinction, no exact frontier. Where does <A> really end and <Non-A> begin? We extend Zadeh's fuzzy set term to fuzzy concept.

**Paradoxist Existentialism**:
life's value consists in its lack of value;
life's sense consists in its lack of sense.

**Semantic Paradox (I)**: I AM WHO I AM NOT.
If I am not Socrates, and since I am who I am not, it results that I am Socrates.
If I am Socrates, and since I am who I am not, it results that I am not Socrates.
Generally speaking: "I am X" if and only if "I am not X".
Who am I?
In a similar pattern one constructs the paradoxes:
I AM MYSELF WHEN I AM NOT MYSELF.
I EXIST WHEN I DON'T EXIST.
And, for the most part:
I {verb} WHEN I DON'T {verb}.
(F. Smarandache, "Linguistic Paradoxes")

What is a dogma?
An idea that makes you have no other idea.
How can we get rid of such authoritative tenet? [To un-read and un-study it!]

**Semantic Paradox (II):** I DON'T THINK.
This can not be true for, in order to even write this sentence, I needed to think (otherwise I was writing with mistakes, or was not writing it at all).
Whence "I don't think" is false, which means "I think".

**Unsolved Mysteries**:
a) Is it true that for each question there is at least an answer?
b) Is any statement the result of a question?
c) Let P(n) be the following assertion:
"If S(n) is true, then S(n+1) is false", where S(n) is a sentence relating on parameter n.
Can we prove by mathematical induction that P(n) is true?
d) "<A> is true if and only if <A> is false".
Is this true or false?
e) How can this assertion "Living without living" be true?
Find a context. Explain.

<Anti-A> of <A>.
Anti-literature of literature.
<Non-A> of <A>.
Language of non-language.

<A> of <Non-A>.
Artistic of the non-artistic.

**Tautologies**:
I want because I want. (showing will, ambition)
<A> because of <A>.
(F. Smarandache, "Linguistic Tautologies")

Our axiom is to break down all axioms.

Be patient without patience.

The non-existence exists.
The culture exists through its non-existence.
Our culture is our lack of culture.

Style without style.

The rule we apply:  it is no rule.

**Paradox of the Paradoxes**:
Is "This is a paradox" a paradox?
I mean is it true or false?

To speak without speaking.  Without words  (body language).
To communicate without communicating.
To do the un-do.

To know nothing about everything, and everything about nothing.

I do only what I can't!
If I can't do something, of course "I can do" is false.
And, if I can do, it's also false because I can do only what I am not able to do.

I cannot for I can.

Paradoxal sleep, from a French "Larousse" dictionary (1989), is a phase of the sleep when the dreams occur.
Sleep, sleep, but why paradoxal?
How do the dreams put up with reality?

Is O. J. Simpson's crime trial an example of:  justice of injustice, or injustice of justice?
However, his famous release is a victory against the system!

Corrupt the incorruptible!

Everything, which is not paradoxist, is however paradoxist.  This is the Great Universal Paradox.
A superparadox;
(as a superman in a hyperspace).

Facts exist in isolation from other facts (= the analytic philosophy),
and in connection as well with each other (= Whitehead's and Bergson's thoughts).
The neutrosophic philosophy unifies contradictory and noncontradictory ideas in any human field.

The antagonism doesn't exist.
Or, if the antagonism does exist, this becomes (by neutrosophic view) a non-(or un-)antagonism:  a normal thought.  I don't worry about it as well as Wordsworth.

Platonism is the observable of unobservable, the thought of the non-thought.

The essence of a thing may never be reached.  It is a symbol, a pure and abstract and absolute notion.

An action may be considered g% good (or right) and b% bad (or wrong), where $0 \le g+b \le 100$ - the remainder being indeterminacy, not only <good> or only <bad> - with rare exceptions, if its consequence is g% happiness (pleasure).
In this case the action is g%-useful (in a *semi-utilitarian* way).
Utilitarism shouldn't work with absolute values only!

Verification has a pluri-sense because we have to demonstrate or prove that something is t% true, and f% false, where $0 \le t,f \le 100$ and $t+f \le 100$, not only t = 0 or 100 - which occurs in rare/absolute exceptions, by means of formal rules of reasoning of this neutrosophic philosophy.

The logical cogitation's structure is discordant.
Scientism and Empiricism are strongly related. They can't run one without other, because one exists in order to complement the other and to differentiate it from its opponent.
PLUS doesn't work without MINUS, and both of them supported by ZERO. They all are cross-penetrating sometimes up to confusion.
The non-understandable is understandable.
If vices wouldn't exist, the virtues will not be seen (T. Muşatescu).
Any new born theory (notion, term, event, phenomenon) automatically generates its non-theory - not necessarily anti-(notion, term, event, phenomenon). Generally speaking, for any <A> a <Non-A> (not necessarily <Anti-A>) will exist for compensation.
The neutrosophy is a theory of theories, because at any moment new ideas and conceptions are appearing and implicitly their negative and neutral senses are highlighted.
Connections & InterConnections...

The non-important is important, because the first one is second one's shadow that makes it grow its value.
The important things would not be so without any unimportant comparison.

The neutrosophic philosophy accepts *a priori* & *a posteriori* any philosophical idea, but associates it with adverse and neutral ones, as a *summum*.
This is to be neutrosophic without being!
Its schemes are related to the neutrality of everything.

## III) ON RUGINA'S ORIENTATION TABLE

Starting from a new viewpoint in philosophy, the neutrosophy, one extends the classical 'probability theory', 'fuzzy set' and 'fuzzy logic' to <neutrosophic probability>, <neutrosophic set> and <neutrosophic logic> respectively.

They are useful in artificial intelligence, neural networks, evolutionary programming, neutrosophic dynamic systems, quantum theory, and decision making in economics.
With the neutrosophic logic help one explores Rugina's Orientation Table, a remarkable tool of study, at the micro- and macro-level, of problems in all sciences.

### 1) NEUTROSOPHY, A NEW BRANCH OF MATHEMATICAL PHILOSOPHY

A) <u>Etymology</u>:
**Neutro-sophy** [French *neutre* < Latin *neuter*, neutral, and Greek *sophia*, skill/wisdom] means knowledge of neutral thought.

B) <u>Definition</u>:
**Neutrosophy** is a new branch of philosophy which studies the origin, nature, and scope of neutralities, as well as their interactions with different ideational spectra.

**C) Characteristics**:

This mode of thinking:

- proposes new philosophical theses, principles, laws, methods, formulas, movements;
- interprets the uninterpretable;
- regards, from many different angles, old concepts, systems:

showing that an idea, which is true in a given referential system, may be false in another one, and vice versa;

- measures the stability of unstable systems,

and instability of stable systems.

**D) Methods of Neutrosophic Study**:

mathematization (neutrosophic logic, neutrosophic probability and statistics, duality), generalization, complementarity, contradiction, paradox, tautology, analogy, reinterpretation, combination, interference, aphoristic, linguistic, multidisciplinarity.

**E) Formalization**:

Let's note by <A> an idea or theory or concept, by <Non-A> what is not <A>, and by <Anti-A> the opposite of <A>. Also, <Neut-A> means what is neither <A>, nor <Anti-A>, i.e. neutrality in between the two extremes. And <A'> a version of <A>.

<Non-A> is different from <Anti-A>.

For example:

If <A> = white, then <Anti-A> = black (antonym),

but <Non-A> = green, red, blue, yellow, black, etc. (any color, except white), while <Neut-A> = green, red, blue, yellow, etc. (any color, except white and black), and <A'> = dark white, etc. (any shade of white).

<Neut-A> ≡ <Neut-(Anti-A)>, neutralities of <A> are identical with neutralities of <Anti-A>.

<Non-A> ⊃ <Anti-A>, and <Non-A> ⊃ <Neut-A> as well,

also

<A> ∩ <Anti-A> = ∅,

<A> ∩ <Non-A> = ∅.

<A>, <Neut-A>, and <Anti-A> are disjoint two by two.

<Non-A> is the completitude of <A> with respect to the universal set.

**F) Main Principle**:

Between an idea <A> and its opposite <Anti-A>, there is a continuum-power spectrum of neutralities <Neut-A>.

**G) Fundamental Thesis**:

Any idea <A> is t% true, i% indeterminate, and f% false,

where t+i+f = 100.

**H) Main Laws**:

Let <Æ> be an attribute, and (a, i, b) ∈ $[0, 100]^3$,

with a+i+b = 100. Then:

- There is a proposition <P> and a referential system <*R*>,

such that <P> is a% <Æ>, i% indeterminate or <Neut-Æ>, and b% <Anti-Æ>.

- For any proposition <P>, there is a referential system <*R*>, such that <P> is a% <Æ>, i% indeterminate or <Neut-Æ>, and b% <Anti-Æ>.
- <Æ> is at some degree <Anti-Æ>, while <Anti-Æ> is at some degree <Æ>.

## 2) NEUTROSOPHIC PROBABILITY AND NEUTROSOPHIC STATISTICS

Let's first generalize the classical notions of 'probability' and 'statistics' for practical reasons.

A) Definitions:

**Neutrosophic Probability** studies the chance that a particular event E will occur, where that chance is represented by three coordinates (variables): t% true, i% indeterminate, and f% false, with $t+i+f = 100$ and $f,i,t \in [0, 100]$.

**Neutrosophic Statistics** is the analysis of such events.

B) Neutrosophic Probability Space:

The universal set, endowed with a neutrosophic probability defined for each of its subset, forms a neutrosophic probability space.

C) Applications:

1) The probability that candidate C will win an election is say 25% true (percent of people voting for him), 35% false (percent of people voting against him), and 40% indeterminate (percent of people not coming to the ballot box, or giving a blank vote - not selecting anyone, or giving a negative vote - cutting all candidates on the list).

Dialectic and dualism don't work in this case anymore.

2) Another example, the probability that tomorrow it will rain is say 50% true according to meteorologists who have investigated the past years' weather, 30% false according to today's very sunny and droughty summer, and 20% undecided (indeterminate).

## 3) NEUTROSOPHIC SET

Let's second generalize, in the same way, the fuzzy set.

A) Definition:

**Neutrosophic Set** is a set such that an element belongs to the set with a neutrosophic probability, i.e. t% is true that the element is in the set, f% false, and i% indeterminate.

B) Neutrosophic Set Operations:

Let M and N be two neutrosophic sets.

One can say, by language abuse, that any element neutrosophically belongs to any set, due to the percentage of truth/indeterminacy/falsity which varies between 0 and 100.

For example: $x(50,20,30) \in M$ (which means, with a probability of 50% x is in M, with a probability of 30% x is not in M, and the rest is undecidable), or $y(0,0,100) \in M$ (which normally means y is not for sure in M), or $z(0,100,0) \in M$ (which means one doesn't know absolutely anything about z's affiliation with M).

Let $0 \leq t_1, t_2, t' \leq 1$ represent the truth-probabilities,
$0 \leq i_1, i_2, i' \leq 1$ the indeterminacy-probabilities, and
$0 \leq f_1, f_2, f' \leq 1$ the falsity-probabilities of an element x to be in the set M and in the set N respectively, and of an element y to be in the set N, where $t_1 + i_1 + f_1 = 1$, $t_2 + i_2 + f_2 = 1$, and $t' + i' + f' = 1$.

One notes, with respect to the given sets,

$x = x(t_1, i_1, f_1) \in M$ and $x = x(t_2, i_2, f_2) \in N$,

by mentioning x's neutrosophic probability appurtenance.

And, similarly, $y = y(t', i', f') \in N$.

Also, for any $0 \leq x \leq 1$ one notes $1-x = \bar{x}$.

Let $W(a,b,c) = (1-a)/(b+c)$ and $W(R) = W(R(t),R(i),R(f))$ for any tridimensional vector $R = ( R(t),R(i),R(f) )$.

Complement of M:

Let $N(x) = 1-x = \bar{x}$. Therefore:

if $x( t_1, i_1, f_1 ) \in M$,

then $x( N(t_1), N(i_1)W(N), N(f_1)W(N) ) \in C(M)$.

Intersection:

Let $C(x,y) = xy$, and $C(z_1,z_2) = C(z)$ for any bidimensional vector $z = (z_1, z_2)$. Therefore:
if $x( t_1, i_1, f_1 ) \in M$, $x( t_2, i_2, f_2 ) \in N$,
then $x( C(t), C(i)W(C), C(f)W(C) ) \in M \cap N$.

Union:

Let $D1(x,y) = x+y-xy = x+\bar{x}y = y+x\bar{y}$, and $D1(z_1,z_2) = D1(z)$ for any bidimensional vector $z = (z_1, z_2)$. Therefore:
if $x( t_1, i_1, f_1 ) \in M$, $x( t_2, i_2, f_2 ) \in N$,
then $x( D1(t), D1(i)W(D1), D1(f)W(D1) ) \in M \cup N$.

Cartesian Product:

if $x( t_1, i_1, f_1 ) \in M$, $y( t', i', f' ) \in N$,
then $( x( t_1, i_1, f_1 ), y( t', i', f' ) ) \in M \times N$.

Difference:

Let $D(x,y) = x-xy = x\bar{y}$, and $D(z_1,z_2) = D(z)$ for any bidimensional vector $z = (z_1, z_2)$. Therefore:
if $x( t_1, i_1, f_1 ) \in M$, $x( t_2, i_2, f_2 ) \in N$,
then $x( D(t), D(i)W(D), D(f)W(D) ) \in M \setminus N$,
because $M \setminus N = M \cap C(N)$.

C) Applications:

From a pool of refugees, waiting in a political refugee camp to get the America visa of emigration, a% are accepted, r% rejected, and p% in pending (not yet decided), a+r+p=100. The chance of someone in the pool to emigrate to USA is not a% as in classical probability, but a% true and p% pending (therefore normally bigger than a%) - because later, the p% pending refugees will be distributed into the first two categories, either accepted or rejected.

Another example, a cloud is a neutrosophic set, because its borders are ambiguous, and each element (water drop) belongs with a neutrosophic probability to the set (e.g. there are separated water drops, around a compact mass of water drops, that we don't know how to consider them: in or out of the cloud).
We are not sure where the cloud ends nor where it begins, neither if some elements are or are not in the set. That's why the percent of indeterminacy is required: for a more organic, smooth, and especially accurate estimation.

## 4) NEUTROSOPHIC LOGIC, A GENERALIZATION OF FUZZY LOGIC

A) Introduction:

One passes from the classical {0, 1} Bivalent Logic of George Boole, to the Three-Valued Logic of Reichenbach (leader of the logical empiricism), then to the $\{0, a_1, ..., a_n, 1\}$ Plurivalent one of Łukasiewicz (and Post's m-valued calculus), and finally to the [0, 1] Infinite Logic as in mathematical analysis and probability: a Transcendental Logic (with values of the power of continuum), or Fuzzy Logic.

Falsehood is infinite, and truthhood quite alike; in between, at different degrees, indeterminacy as well.
Everything is G% good, I% indeterminate, and B% bad,
where G + I + B = 100.

Besides Diderot's dialectics on good and bad ("Rameau's Nephew", 1772), any act has its percentage of "good", "indeterminate", and of "bad" as well incorporated.

Rodolph Carnap said:
"Metaphysical propositions are neither true nor false, because they assert nothing, they contain neither knowledge nor error (...)".
Hence, there are infinitely many statuses in between "Good" and "Bad", and generally speaking in between "A" and "Anti-A", like on the real number segment:

| | |
|---|---|
| [0, | 1] |
| False | True |
| Bad | Good |
| Non-sense | Sense |
| Anti-A | A |

0 is the absolute falsity, 1 the absolute truth. In between each opposing pair, normally in a vicinity of 0.5, are being set up the neutralities.

There exist as many states in between "True" and "False" as in between "Good" and "Bad". Irrational and transcendental standpoints belong to this interval.

Even if an act apparently looks to be only good, or only bad, the other hidden side should be sought. The ratios

$$\frac{\text{Anti-A}}{A}, \quad \frac{\text{Non-A}}{A}$$

vary indefinitely. They are transfinite.

If a statement is 30%T (true) and 60%I (indeterminate), then it is 10%F (false). This is somehow alethic, meaning pertaining to truthhood and falsehood in the same time.

In opposition to Fuzzy Logic, if a statement is 30%T doesn't involve it is 70%F. We have to study its indeterminacy as well.

B) <u>Definition of Neutrosophic Logic</u>:

This is a generalization (for the case of null indeterminacy) of the fuzzy logic.

Neutrosophic logic is useful in the real-world systems for designing control logic, and may work in quantum mechanics.

If a proposition P is t% true, doesn't necessarily mean it is 100-t% false as in fuzzy logic. There should also be a percent of indeterminacy on the values of P.

A better approach of the logical value of P is f% false, i% indeterminate, and t% true, where $t+i+f = 100$ and $t,i,f \in [0, 100]$, called neutrosophic logical value of P, and noted by $n(P) = (t,i,f)$.

**Neutrosophic Logic** means the study of neutrosophic logical values of the propositions.

There exist, for each individual event, PRO parameters, CONTRA parameters, and NEUTER parameters which influence the above values.

Indeterminacy results from any hazard which may occur, from unknown parameters, or from new arising conditions.

This resulted from practice.

C) <u>Applications</u>:

1) The candidate C, who runs for election in a metropolis M of p people with right to vote, will win.

This proposition is, say, 25% true (percent of people voting for him), 35% false (percent of people voting against him), and 40% indeterminate (percent of people not coming to the ballot box, or giving a blank vote - not selecting anyone, or giving a negative vote - cutting all candidates on the list).

2) Tomorrow it will rain.

This proposition is, say, 50% true according to meteorologists who have investigated the past years' weather, 30% false according to today's very sunny and droughty summer, and 20% undecided.

3) This is a heap.

As an application to the sorites paradoxes, we may now say this proposition is t% true, f% false, and i% indeterminate (the neutrality comes for we don't know exactly where is the difference between a heap and a non-heap; and, if we approximate the border, our 'accuracy' is subjective).

We are not able to distinguish the difference between yellow and red as well if a continuum spectrum of colors is painted on a wall imperceptibly changing from one into another.

D) <u>Definition of Neutrosophic Logical Connectors</u>:

One uses the definitions of neutrosophic probability and neutrosophic set.

Let $0 \le t_1, t_2 \le 1$ represent the truth-probabilities,

$0 \le i_1, i_2 \le 1$ the indeterminacy-probabilities, and

$0 \le f_1, f_2 \le 1$ the falsity-probabilities of two events $P_1$ and $P_2$ respectively, where $t_1 + i_1 + f_1 = 1$ and $t_2 + i_2 + f_2 = 1$.

One notes the neutrosophic logical values of $P_1$ and $P_2$ by

$n(P_1) = (t_1, i_1, f_1)$ and $n(P_2) = (t_2, i_2, f_2)$.

Also, for any $0 \le x \le 1$ one notes $1\text{-}x = \bar{x}$.

Let $W(a,b,c) = (1\text{-}a)/(b+c)$ and $W(R) = W(R(t),R(i),R(f))$ for any tridimensional vector $R = ( R(t),R(i),R(f) )$.

Negation:

Let $N(x) = 1\text{-}x = \bar{x}$. Then:

$$n(\neg P_1) = ( N(t_1), N(i_1)W(N), N(f_1)W(N) ).$$

Conjunction:

Let $C(x,y) = xy$, and $C(z_1,z_2) = C(z)$ for any bidimensional vector $z = (z_1, z_2)$. Then:

$$n(P_1 \wedge P_2) = ( C(t), C(i)W(C), C(f)W(C) ).$$

(And, in a similar way, generalized for n propositions.)

Weak or inclusive disjunction:

Let $D1(x,y) = x+y\text{-}xy = x+\bar{x}y = y+x\bar{y}$, and $D1(z_1,z_2) = D1(z)$ for any bidimensional vector $z = (z_1, z_2)$. Then:

$$n(P_1 \vee P_2) = ( D1(t), D1(i)W(D1), D1(f)W(D1) ).$$

(And, in a similar way, generalized for n propositions.)

Strong or exclusive disjunction:

Let $D2(x,y) = x(1\text{-}y)+y(1\text{-}x)\text{-}xy(1\text{-}x)(1\text{-}y) = x\bar{y}+\bar{x}y\text{-}xy\overline{xy}$, and $D2(z_1,z_2) = D2(z)$ for any bidimensional vector $z = (z_1, z_2)$. Then:

$$n(P_1 \underline{\vee} P_2) = ( D2(t), D2(i)W(D2), D2(f)W(D2) ).$$

(And, in a similar way, generalized for n propositions.)

Material conditional (implication):

Let $I(x,y) = 1\text{-}x+xy = \bar{x}+xy = 1\text{-}x\bar{y}$, and $I(z1,z2) = I(z)$ for any bidimensional vector $z = (z_1, z_2)$. Then:

$$n(P_1 \rightarrow P_2) = ( I(t), I(i)W(I), I(f)W(I) ).$$

Material biconditional (equivalence):

Let $E(x,y) = (1\text{-}x+xy)(1\text{-}y+xy) = (\bar{x}+xy)(\bar{y}+xy) = (1\text{-}x\bar{y})(1\text{-}\bar{x}y)$, and

$E(z_1,z_2) = E(z)$ for any bidimensional vector $z = (z_1, z_2)$.

$$n(P \leftrightarrow Q) = ( E(t), E(i)W(E), E(f)W(E) ).$$

Sheffer's connector:

Let $S(x,y) = 1\text{-}xy$, and $S(z_1,z_2) = S(z)$ for any bidimensional vector $z = (z_1, z_2)$.

$$n(P|Q) = n(\neg P \vee \neg Q) = ( S(t), S(i)W(S), S(f)W(S) ).$$

Peirce's connector:

Let $P(x,y) = (1\text{-}x)(1\text{-}y) = \overline{xy}$, and $P(z_1,z_2) = P(z)$ for any bidimensional vector $z = (z_1, z_2)$.

$$n(P \downarrow Q) = n(\neg P \wedge \neg Q) = ( P(t), P(i)W(P), P(f)W(P) ).$$

E) <u>Properties of Neutrosophic Logical Connectors</u>:

Let's note by t(P) the truth-component of the neutrosophic value n(P), and t(P) = p, t(Q) = q.

a) Conjunction:

$t(P \wedge Q) \leq \min \{p, q\}$.

$$\bigwedge_{k=1}^{\infty} t(P) = 0 \text{ if } t(P) \neq 1.$$

b) Weak disjunction:

$t(P \vee Q) \geq \max \{p, q\}$.

$$\bigvee_{k=1}^{\infty} t(P) = 1 \text{ if } t(P) \neq 0.$$

c) Implication:

$t(P \to P) = 1$ if $t(P) = 0$ or 1, and $\rangle$ p otherwise.

$$\lim_{t(P)\to 0} t(P \to Q) = 1$$

$$\lim_{t(Q)\to 1} t(P \to Q) = 1$$

$$\lim_{t(P)\to 1} t(P \to Q) = q$$

$$\lim_{t(Q)\to 0} t(P \to Q) = 1-p$$

d) Equivalence:

$t(P \leftrightarrow Q) = t(Q \leftrightarrow P) = t(\neg P \leftrightarrow \neg Q)$

$$\lim_{\substack{t(P)\to 0\\ t(Q)\to 0}} t(P \leftrightarrow Q) = 1$$

$$\lim_{\substack{t(P)\to 1\\ t(Q)\to 1}} t(P \leftrightarrow Q) = 1$$

$$\lim_{\substack{t(P)\to 0\\ t(Q)\to 1}} t(P \leftrightarrow Q) = 0$$

$$\lim_{t(P)\to 1} t(P \leftrightarrow Q) = 0$$

$t(Q)\to 0$

$\lim_{t(P)\to 0} t(P \leftrightarrow Q) = 1-q$

$\lim_{t(P)\to 1} t(P \leftrightarrow Q) = q$

Let $q \neq 0,1$ be constant, and one notes
$p_{max}(q) = (q^2-3q+1)/(2q^2-2q)$. Then:

$\max_{0 \leq t(P) \leq 1} t(P \leftrightarrow Q)$ occurs when:

$p = p_{max}(q)$ if $p_{max}(q) \in [0, 1]$,
or $p = 0$ if $p_{max}(q) < 0$,
or $p = 1$ if $p_{max}(q) > 1$,
because the equivalence connector is described by a parabola of equation

$e_q(p) = (q^2-q)p^2 + (-q^2+3q-1)p + (1-q)$,

which is concave down.

## 5) NEUTROSOPHIC TOPOLOGY

A) <u>Definition</u>:

Let's construct a **Neutrosophic Topology** on NT = [0, 1], considering the associated family of subsets (0, p), for $0 \leq p \leq 1$, the whole set [0, 1], and the empty set $\emptyset = (0, 0)$, called open sets, which is closed under set union and finite intersection. The union is defined as $(0, p) \cup (0, q) = (0, d)$, where $d = p+q-pq$, and the intersection as $(0, p) \cap (0, q) = (0, c)$, where $c = pq$.
The complementary of (0, p) is (0, n), where $n = 1-p$, which is a closed set.

B) <u>Neutrosophic Topological Space</u>:

The interval NT, endowed with this topology, forms a neutrosophic topological space.

C) <u>Isomorphicity</u>:

Neutrosophic Logical Space, Neutrosophic Topological Space, and Neutrosophic Probability Space are all isomorphic.

A method of Neutrosophy is the:

## 6) TRANSDISCIPLINARITY:

**A) <u>Introduction</u>**:

Transdisciplinarity means to find common features to uncommon entities: $\langle A \rangle \cap \langle Non\text{-}A \rangle \neq \emptyset$, even if they are disjunct.

**B) <u>Multi-Structure and Multi-Space</u>**:

I consider that life and practice do not deal with 'pure' spaces, but with a group of many spaces,
with a mixture of structures, a 'mongrel', a heterogeneity --
the ardently preoccupation is to reunite them, to constitute a
multi-structure.

I thought to a multi-space also: fragments (potsherds) of spaces put together, say as an example: Banach, Hausdorff, Tikhonov, compact, paracompact, Fock symmetric, Fock antisymmetric, path-connected, simply connected, discrete metric, indiscrete pseudo-metric, etc. spaces that work together as a whole mechanism. The difficulty is to be the passage over 'frontiers' (borders between two disjoint spaces);

i.e. how can we organically tie a point P1 from a space S1 with
a point P2 from a structurally opposite space S2 ?
Does the problem become more complicated when the spaces' sets are not disjoint?

Let $S_1$ and $S_2$ be two distinct structures, induced by the group of laws L which verify the axiom groups $A_1$ and $A_2$ respectively, such that $A_1$ is strictly included in $A_2$.
One says that the set M, endowed with the properties:
a) M has an $S_1$-structure,
b) there is a proper subset P (different from the empty set, from the unitary element, and from M) of the initial set M which has an $S_2$-structure,
c) M doesn't have an $S_2$-structure,
is called an **$S_1$-structure with respect to the $S_2$-structure**.

Let $S_1, S_2, ..., S_k$ be distinct space-structures.
We define the **Multi-Space** (or **k-structured-space**) as a set M such that for each structure $S_i$, $1 \leq i \leq k$, there is a proper (different from $\emptyset$ and from M) subset $M_i$ of it which has that structure. The $M_1, M_2, ..., M_k$ proper subsets are different two by two.

Let's introduce new terms:

**C) <u>Psychomathematics</u>**:
A discipline which studies psychological processes in connection with mathematics.

**D) <u>Mathematical Modeling of Psychological Processes</u>**:
Weber's law and Fechner's law on sensations and stimuli are improved.

**E) <u>Psychoneutrosophy</u>**:
Psychology of neutral thought, action, behavior, sensation, perception, etc. This is a hybrid field deriving from theology, philosophy, economics, psychology, etc.
For example, to find the psychological causes and effects of individuals supporting neutral ideologies (neither capitalists, nor communists), politics (not in the left, not in the right), etc.

**F) <u>Socioneutrosophy</u>**:
Sociology of neutralities.
For example the sociological phenomena and reasons which determine a country or group of people or class to remain neuter in a military, political, ideological, cultural, artistic, scientific, economical, etc. international or internal war (dispute).

**G) <u>Econoneutrosophy</u>**:
Economics of non-profit organizations, groups, such as: churches, philanthropic associations, charities, emigrating foundations, artistic or scientific societies, etc.
How they function, how they survive, who benefits and who loses, why are they necessary, how they improve, how they interact with for-profit companies.

These terms are in the process of development.

## 6) RUGINA'S ORIENTATION TABLE

In order to clarify the anomalies in science, Rugina (1989, 1998) proposes an original method, starting first from an economic point of view but generalizing it to any science, to study the equilibrium and disequilibrium of systems. His Table comprises seven basic models:

Model $M_1$ (which is 100% stable),
Model $M_2$ (which is 95% stable, and 5% unstable),
Model $M_3$ (which is 65% stable, and 35% unstable),
Model $M_4$ (which is 50% stable, and 50% unstable),
Model $M_5$ (which is 35% stable, and 65% unstable),
Model $M_6$ (which is 100% unstable).

He gives Orientation Tables for Physical Sciences and Mechanics (Rugina, 1989, p. 18), for the Theory of Probability, for Logic, and generally for any Natural or Social Science (Rugina, 1989, pp. 286-288).

"An anomaly can be simply defined as a deviation from a position of stable equilibrium represented by Model $M_1$" (Rugina, 1989, p. 17).

Rugina proposes the Universal Hypothesis of Duality:

"The physical universe in which we are living, including human society and the world of ideas, all are composed in different and changeable proportions of stable (equilibrium) and unstable (disequilibrium) elements, forces, institutions, behavior and value"

and the General Possibility Theorem:

"there is an unlimited number of possible combinations or systems in logic and other sciences".

According to the last assertations one can extend Rugina's Orientation Table in the way that any system in each science is s% stable and u% unstable, with s+u=100 and both parameters
$0 \leq s, u \leq 100$, somehow getting to a fuzzy approach.

But, because each system has hidden features and behaviors, and there would always be unexpected occuring conditions we are not able to control - we mean the indeterminacy plays a role as well, a better approach would be the *Neutrosophic Model*:

Any system in each science is s% stable, i% indeterminate, and u% unstable, with s+i+u=100 and all three parameters $0 \leq s,i,u \leq 100$.

EXAMPLE OF MODEL M3 IN RUGINA'S ORIENTATION TABLE:

<u>The Paradoxist Geometry</u> (actually the percentage of instability is between 20-35)

In 1969, intrigued by geometry, I simultaneously constructed a partially Euclidean and partially Non-Euclidean space by a strange replacement of the Euclid's fifth postulate (axiom of parallels) with the following five-statement proposition:

This geometry unites all together: Euclid, Lobachevsky/Bolyai, and Riemann Geometries. And separates them as well!

a) there are at least a straight line and a point exterior
to it in this space for which only one line passes through the point
and does not intersect the initial line;
[1 parallel]

b) there are at least a straight line and a point exterior
to it in this space for which only a finite number of
lines $l_1, ..., l_k$ ($k >= 2$) passe through the point and do not
intersect the initial line;
[2 or more (in a finite number) parallels]

c) there are at least a straight line and a point exterior
to it in this space for which any line that passes through the point
intersects the initial line;
[0 parallels]

d) there are at least a straight line and a point exterior
to it in this space for which an infinite number of lines
that passe through the point (but not all of them) do not intersect
the initial line;

[an infinite number of parallels, but not all lines passing through]

e) there are at least a straight line and a point exterior to it in this space for which any line that passes through the point does not intersect the initial line;
[an infinite number of parallels, all lines passing through the point].

FIRST EXAMPLE OF MODEL M7 IN RUGINA'S ORIENTATION TABLE:
**The Non-Geometry** (the percentage of instability is 100)

It's a lot easier to deny the Euclid's five postulates than Hilbert's twenty thorough axioms.

1. It is not always possible to draw a line from an arbitrary point to another arbitrary point.

   For example:
   this axiom can be denied only if the model's space has at least a discontinuity point;
   (in our bellow model MD, one takes an isolated point I in between f1 and f2, the only one which will not verify the axiom).

2. It is not always possible to extend by continuity a finite line to an infinite line.

   For example:
   consider the bellow Model, and the segment AB, where both A and B lie on f1, A in between P and N, while B on the left side of N;
   one can not at all extend AB either beyond A or beyond B, because the resulted curve, noted say A'-A-B-B', would not be a geodesic (i.e. line in our Model) anymore.

   If A and B lie in delta1-f1, both of them closer to f1, A in the left side of P, while B in the right side of P,
   then the segment AB, which is in fact A-P-B, can be extended beyond A and also beyond B only up to f1
   (therefore one gets a finite line too, A'-A-P-B-B', where A', B' are the intersections of PA, PB respectively with f1).

   If A, B lie in delta1-f1, far enough from f1 and P, such that AB is parallel to f1, then AB verifies this postulate.

3. It is not always possible to draw a circle from an arbitrary point and of an arbitrary interval.

   For example:
   same as for the first axiom;
   the isolated point I, and a very small interval not reaching f1 neither f2, will deny this axiom.

4. Not all the right angles are congruent.

   (See example of the Anti-Geometry, explained bellow.)

5. If a line, cutting two other lines, forms the interior angles of the same side of it strictly less than two right angles, then not always the two lines extended towards infinite cut each other in the side where the angles are strictly less than two right angles.

For example:
let h1, h2, l be three lines in delta1-delta2, where h1 intersects f1 in A, and h2 intersects f1 in B, with A, B, P different each other, such that h1 and h2 do not intersect, but l cuts h1 and h2 and forms the interior angles of one of its side (towards f1) strictly less than two right angles;
the assumption of the fifth postulate is fulfilled, but the consequence does not hold, because h1 and h2 do not cut each other (they may not be extended beyond A and B respectively, because the lines would not be geodesics anymore).

SECOND EXAMPLE OF MODEL M7 IN RUGINA'S ORIENTATION TABLE:
**The Counter-Projective Geometry** (the percentage of instability is 100)

Let P, L be two sets, and r a relation included in PxL. The elements of P are called points, and those of L lines. When (p, l) belongs to r, we say that the line l contains the point p.
For these, one imposes the following COUNTER-AXIOMS:

(a) There exist: either at least two lines, or no line, that contains two given distinct points.

(b) Let p1, p2, p3 be three non-collinear points, and q1, q2 two distinct points. Suppose that {p1, q1, p3} and {p2, q2, p3} are collinear triples. Then the line containing p1, p2, and the line containing q1, q2 do not intersect.

(c) Every line contains at most two distinct points.

Does the Duality Principle hold in a counter-projective space?
What about Desargues's Theorem, Fundamental Theorem of Projective Geometry / Theorem of Pappus, and Staudt Algebra ?
Or Pascal's Theorem, Brianchon's Theorem ?
(I think none of them will hold!)
However Rugina's Hypothesis of Duality does hold (although the this geometry is formed by unstable elements only!).

The Theory of Buildings of Tits, which contains the Projective Geometry as a particular case, can be 'distorted' in the same <paradoxist> way by deforming its axiom of a BN-pair (or Tits system) for the triple (G, B, N), where G is a group, and B, N its subgroups;
[see J. Tits, "Buildings of spherical type and finite BN-pairs", Lecture Notes in Math 386, Springer, 1974].
Notions as: simplex, complex, chamber, codimension, apartment, building will get contorted either...
Develop a Theory of Distorted Buildings of Tits!

THIRD EXAMPLE OF MODEL M7 IN RUGINA'S ORIENTATION TABLE:
**The Anti-Geometry** (the percentage of instability is 100 - even... more, this is the geometry of total chaos!)

It is possible to entirely de-formalize Hilbert's groups of axioms of the Euclidean Geometry, and to construct a model such that none of his fixed axioms holds.

Let's consider the following things:

- a set of <points>: A, B, C, ...
- a set of <lines>: h, k, l, ...
- a set of <planes>: alpha, beta, gamma, ...

and

- a set of relationships among these elements: "are situated", "between", "parallel", "congruent", "continuous", etc.

Then, we can deny all Hilbert's twenty axioms [see David Hilbert, "Foundations of Geometry", translated by E. J. Townsend, 1950; and Roberto Bonola, "Non-Euclidean Geometry", 1938].
There exist cases, within a geometric model, when the same axiom is verified by certain points/lines/planes and denied by others.

GROUP I. ANTI-AXIOMS OF CONNECTION:

I.1. Two distinct points A and B do not always completely determine a line.

Let's consider the following model MD:
get an ordinary plane delta, but with an infinite hole in of the following shape:

```
                         P
p                        .
     semi-plane delta1 . .
l                       .  .
                       .    .     curve f1 (frontier)
a                     .      .
   ..........................   ........................
                      N
n
                              .I
e                                .J


                               .K
d

e

l
    ..........................   ..........................
t                     .    .      curve f2 (frontier)
                       .   .
a    semi-plane delta2 .  .
                          .
                         Q
```

Plane delta is a reunion of two disjoint planar semi-planes;
f1 lies in MD, but f2 does not;
P, Q are two extreme points on f that belong to MD.

One defines a LINE l as a geodesic curve: if two points A, B that belong to MD lie in l, then the shortest curve lied in MD between A and B lies in l also.
If a line passes two times through the same point, then it is called double point (KNOT).

One defines a PLANE alpha as a surface such that for any two points A, B that lie in alpha and belong to MD there is a geodesic which passes through A, B and lies in alpha also.

Now, let's have two strings of the same length:
one ties P and Q with the first string s1 such that the curve s1 is folded in two or more different planes and s1 is under the plane delta;
next, do the same with string s2, tie Q with P, but over the plane delta and such that s2 has a different form from s1;
and a third string s3, from P to Q, much longer than s1.
s1, s2, s3 belong to MD.

Let I, J, K be three isolated points -- as some islands, i.e. not joined with any other point of MD, exterior to the plane delta.

This model has a measure, because the (pseudo-)line is the shortest way (length) to go from a point to another (when possible).

Of course, this model is not perfect, and is far from the best. Readers are asked to improve it, or to make up a new one that is better.

(Let A, B be two distinct points in delta1-f1. P and Q are two points on s1, but they do not completely determine a line, referring to the first axiom of Hilbert, because A-P-s1-Q are different from B-P-s1-Q.)

I.2. There is at least a line l and at least two distinct points A and B of l, such that A and B do not completely determine the line l.

(Line A-P-s1-Q are not completely determined by P and Q in the previous construction, because CVB-P-s1-Q is another line passing through P and Q too.)

I.3. Three points A, B, C not situated in the same line do not always completely determine a plane alpha.
(Let A, B be two distinct points in delta1-f1, such that A, B, P are not co-linear. There are many planes containing these three points: delta1 extended with any surface s containing s1, but not cutting s2 in between P and Q, for example.)

I.4. There is at least a plane, alpha, and at least three points A, B, C in it not lying in the same line, such that A, B, C do not completely determine the plane alpha.

(See the previous example.)

I.5. If two points A, B of a line l lie in a plane alpha, doesn't mean that every point of l lies in alpha.

(Let A be a point in delta1-f1, and B another point on s1 in between P and Q. Let alpha be the following plane:
delta1 extended with a surface s containing s1, but not cutting s2 in between P and Q, and tangent to delta2 on a line QC, where C is a point in delta2-f2. Let D be point in delta2-f2, not lying on the line QC. Now, A, B, D are lying on the same line A-P-s1-Q-D, A, B are in the plane alpha, but D do not.)

I.6. If two planes alpha, beta have a point A in common, doesn't mean they have at least a second point in common.

(Construct the following plane alpha: a closed surface containing s1 and s2, and intersecting delta1 in one point only, P. Then alpha and delta1 have a single point in common.)

I.7. There exist lines where lies only one point,
or planes where lie only two points,
or space where lie only three points.

(Hilbert's I.7 axiom may be contradicted if the model has discontinuities.
Let's consider the isolated points area.
The point I may be regarded as a line, because it's not possible to add any new point to I to form a line.
One constructs a surface that intersects the model only in the points I and J.)

GROUP II. ANTI-AXIOMS OF ORDER:

II.1. If A, B, C are points of a line and B lies between A and C, doesn't mean that always B lies also between C and A.

[Let T lie in s1, and V lie in s2, both of them closer to Q, but different from it. Then:
P, T, V are points on the line P-s1-Q-s2-P
( i.e. the closed curve that starts from the point P and lies in s1 and passes through the point Q and lies back to s2 and ends in P ),

and T lies between P and V
-- because PT and TV are both geodesics --,
but T doesn't lie between V and P
-- because from V the line goes to P and then to T,
therefore P lies between V and T.]

[By definition: a segment AB is a system of points
lying upon a line between A and B (the extremes are included).

Warning: AB may be different from BA;
for example:
the segment PQ formed by the system of points
starting with P, ending with Q, and lying in s1, is different
from the segment QP formed by the system of points starting
with Q, ending with P,
but belonging to s2.
Worse, AB may be sometimes different from AB;
for example:
the segment PQ formed by the system of points
starting with P, ending with Q, and lying in s1, is different
from the segment PQ formed by the system of points starting
with P, ending with Q, but belonging to s2.]

II.2. If A and C are two points of a line, then:
there does not always exist a point B lying between A and C, or
there does not always exist a point D such that C lies between A
and D.

[For example:
let F be a point on f1, F different from P,
and G a point in delta1, G doesn't belong to f1;
draw the line l which passes through G and F;
then:
there exists a point B lying between G and F
-- because GF is an obvious segment --,
but there is no point D such that F lies between
G and D -- because GF is right bounded in F
( GF may not be extended to the other side of F,
because otherwise the line will not remain a
geodesic anymore ).]

II.3. There exist at least three points situated on
a line such that:
one point lies between the other two,
and another point lies also between the other two.

[For example:
let R, T be two distinct points, different
from P and Q, situated on the line P-s1-Q-s2-P,
such that the lengths PR, RT, TP are all equal;
then:
R lies between P and T,

and T lies between R and P;
also P lies between T and R.]

II.4. Four points A, B, C, D of a line can not always be arranged:
such that B lies between A and C and also between A and D,
and such that C lies between A and D and also between B and D.

[For examples:
- let R, T be two distinct points, different from P and Q, situated on the line P-s1-Q-s2-P such that the lengths PR, RQ, QT, TP are all equal,
therefore R belongs to s1, and T belongs to s2;
then P, R, Q, T are situated on the same line:
such that R lies between P and Q, but not between P and T
-- because the geodesic PT does not pass through R --,
and such that Q does not lie between P and T
-- because the geodesic PT does not pass through Q --,
but lies between R and T;
- let A, B be two points in delta2-f2 such that A, Q, B are collinear, and C, D two points on s1, s2 respectively,
all of the four points being different from P and Q;
then A, B, C, D are points situated on the same line A-Q-s1-P-s2-Q-B, which is the same with line A-Q-s2-P-s1-Q-B, therefore we may have two different orders of these four points in the same time:
A, C, D, B and A, D, C, B.]

II.5. Let A, B, C be three points not lying in the same line, and l a line lying in the same plane ABC and not passing through any of the points A, B, C. Then, if the line l passes through a point of the segment AB, it doesn't mean that always the line l will pass through either a point of the segment BC or a point of the segment AC.

[For example:
let AB be a segment passing through P in the semi-plane delta1, and C a point lying in delta1 too on the left side of the line AB;
thus A, B, C do not lie on the same line;
now, consider the line Q-s2-P-s1-Q-D, where D is a point lying in the semi-plane delta2 not on f2:
therefore this line passes through the point P of the segment AB, but do not pass through any point of the segment BC, nor through any point of the segment AC.]

## GROUP III. ANTI-AXIOM OF PARALLELS.

In a plane alpha there can be drawn through a point A, lying outside of a line l, either no line, or only one line, or a finite number of lines, or an infinite number of lines which do not intersect the line l. (At least two of these situations should occur.)
The line(s) is (are) called the parallel(s) to l through the given point A.

[ For examples:
- let l0 be the line N-P-s1-Q-R, where N is a point lying in delta1 not on f1, and R is a similar point lying in delta2 not on f2, and let A be a point lying on s2, then: no parallel to l0 can be drawn through A (because any line passing through A, hence through s2, will intersect s1, hence l0, in P and Q);
- if the line l1 lies in delta1 such that l1 does not intersect the frontier f1, then: through any point lying on the left side of l1 one and only one parallel will pass;
- let B be a point lying in f1, different from P, and another point C lying in delta1, not on f1; let A be a point lying in delta1 outside of BC; then: an infinite number of parallels to the line BC can be drawn through the point A.

Theorem. There are at least two lines l1, l2 of a plane, which do not meet a third line l3 of the same plane, but they meet each other, ( i.e. if l1 is parallel to l3, and l2 is parallel to l3, and all of them are in the same plane, it's not necessary that l1 is parallel to l2 ).
[ For example:
consider three points A, B, C lying in f1, and different from P, and D a point in delta1 not on f1; draw the lines AD, BE and CE such that E is a point in delta1 not on f1 and both BE and CE do not intersect AD;
then: BE is parallel to AD, CE is also parallel to AD, but BE is not parallel to CE because the point E belong to both of them. ]

## GROUP IV. ANTI-AXIOMS OF CONGRUENCE

IV.1. If A, B are two points on a line l, and A' is a point upon the same or another line l', then:
upon a given side of A' on the line l', we can not always find only one point B' so that the segment AB is congruent to the segment A'B'.

[ For examples:
- let AB be segment lying in delta1 and having no point in common with f1, and construct the line C-P-s1-Q-s2-P (noted by l') which is the same with C-P-s2-Q-s1-P, where C is a point lying in delta1 not on f1 nor on AB;
take a point A' on l', in between C and P, such that A'P is smaller than AB;

now, there exist two distinct points B1' on s1 and B2' on s2, such that A'B1' is congruent to AB and A'B2' is congruent to AB, with A'B1' different from A'B2';
- but if we consider a line l' lying in delta1 and limited by the frontier f1 on the right side (the limit point being noted by M), and take a point A' on l', close to M, such that A'M is less than A'B', then: there is no point B' on the right side of l' so that A'B' is congruent to AB. ]

A segment may not be congruent to itself!

[ For example:
- let A be a point on s1, closer to P, and B a point on s2, closer to P also; A and B are lying on the same line A-Q-B-P-A which is the same with line A-P-B-Q-A, but AB measured on the first representation of the line is strictly greater than AB measured on the second representation of their line. ]

IV.2. If a segment AB is congruent to the segment A'B' and also to the segment A"B", then not always the segment A'B' is congruent to the segment A"B".

[ For example:
- let AB be a segment lying in delta1-f1, and consider the line C-P-s1-Q-s2-P-D, where C, D are two distinct points in delta1-f1 such that C, P, D are collinear. Suppose that the segment AB is congruent to the segment CD (i.e. C-P-s1-Q-s2-P-D).
Get also an obvious segment A'B' in delta1-f1, different from the preceding ones, but congruent to AB.
Then the segment A'B' is not congruent to the segment CD (considered as C-P-D, i.e. not passing through Q).

IV.3. If AB, BC are two segments of the same line l which have no points in common aside from the point B, and A'B', B'C' are two segments of the same line or of another line l' having no point other than B' in common, such that AB is congruent to A'B' and BC is congruent to B'C', then not always the segment AC is congruent to A'C'.

[ For example:
let l be a line lying in delta1, not on f1, and A, B, C three distinct points on l, such that AC is greater than s1;
let l' be the following line: A'-P-s1-Q-s2-P where A' lies in delta1, not on f1, and get B' on s1 such that A'B' is congruent to AB, get C' on s2 such that BC is congruent to B'C' (the points A, B, C are thus chosen);
then: the segment A'C' which is first seen as A'-P-B'-Q-C' is not congruent to AC, because A'C' is the geodesic A'-P-C' (the shortest way from A' to C' does not pass through B') which is strictly less than AC. ]

Definitions. Let h, k be two lines having a point O in common. Then the system (h, O, k) is called the angle of the lines h and k in the point O.

( Because some of our lines are curves, we take the angle of the tangents to the curves in their common point. )

The angle formed by the lines h and k situated in the same plane, noted by <(h, k), is equal to the arithmetic mean of the angles formed by h and k in all their common points.

IV.4. Let an angle (h, k) be given in the plane alpha, and let a line h' be given in the plane beta.
Suppose that in the plane beta a definite side of the line h' be assigned, and a point O'.
Then in the plane beta there are one, or more, or even no half-line(s) k' emanating from the point O' such that the angle (h, k) is congruent to the angle (h', k'), and at the same time the interior points of the angle (h', k') lie upon one or both sides of h'.

[ Examples:

- Let A be a point in delta1-f1, and B, C two distinct points in delta2-f2;
  let h be the line A-P-s1-Q-B, and k be the line A-P-s2-Q-C; because h and k intersect in an infinite number of points (the segment AP), where they normally coincide -- i.e. in each such point their angle is congruent to zero, the angle (h, k) is congruent to zero.
  Now, let A' be a point in delta1-f1, different from A, and B' a point in delta2-f2, different from B, and draw the line h' as A'-P-s1-Q-B';
  there exist an infinite number of lines k', of the form A'-P-s2-Q-C' (where C' is any point in delta2-f2, not on the line QB'), such that the angle (h, k) is congruent to (h', K'), because (h', k') is also congruent to zero, and the line A'-P-s2-Q-C' is different from the line A'-P-s2-Q-D' if D' is not on the line QC'.
- If h, k, and h' are three lines in delta1-P, which intersect the frontier f1 in at most one point, then there exist only one line k' on a given part of h' such that the angle (h, k) is congruent to the angle (h', k').
- *Is there any case when, with these hypotheses, no k' exists ?
- Not every angle is congruent to itself;
  for example:

<(s1, s2) is not congruent to <(s1, s2)
[because one can construct two distinct lines: P-s1-Q-A and P-s2-Q-A, where A is a point in delta2-f2, for the first angle, which becomes equal to zero;
and P-s1-Q-A and P-s2-Q-B, where B is another point in delta2-f2, B different from A, for the second angle, which becomes strictly greater than zero!].

IV. 5. If the angle (h, k) is congruent to the angle (h', k',) and the angle (h'', k''), then the angle (h', k') is not always congruent to the angle (h'', k'').

(A similar construction to the previous one.)

IV. 6. Let ABC and A'B'C' be two triangles such that
AB is congruent to A'B',
AC is congruent to A'C',
<BAC is congruent to <B'A'C'.
Then not always
<ABC is congruent to <A'B'C'
and <ACB is congruent to <A'C'B'.

[For example:
Let M, N be two distinct points in delta2-f2, thus obtaining the triangle PMN;
Now take three points R, M', N' in delta1-f1, such that RM' is congruent to PM, RN' is congruent to RN, and the angle (RM', RN') is congruent to the angle (PM, PN). RM'N' is an obvious triangle.
Of course, the two triangle are not congruent, because for example PM and PN cut each other twice -- in P and Q -- while RM' and RN' only once -- in R.
(These are geodesic triangles.)]

Definitions:

Two angles are called supplementary if they have the same vertex, one side in common, and the other sides not common form a line.

A right angle is an angle congruent to its supplementary angle.

Two triangles are congruent if its angles are congruent two by two, and its sides are congruent two by two.

Propositions:

A right angle is not always congruent to another right angle.

For example:
Let A-P-s1-Q be a line, with A lying in delta1-f1, and B-P-s1-Q another line, with B lying in delta1-f1 and B not lying in the line AP;

we consider the tangent t at s1 in P, and B chosen in a way
that <(AP, t) is not congruent to <(BP, t);
let A', B' be other points lying in delta1-f1
such that <APA' is congruent to <A'P-s1-Q,
and <BPB' is congruent to <B'P-s1-Q.
Then:
- the angle APA' is right, because it is congruent to its supplementary (by construction);
- the angle BPB' is also right, because it is congruent to its supplementary (by construction);
- but <APA' is not congruent to <BPB', because the first one is half of the angle A-P-s1-Q, i.e. half of <(AP, t), while the second one is half of the B-P-s1-Q, i.e. half of <(BP, t).

The theorems of congruence for triangles [side, side, and angle in between;  angle, angle, and common side;  side, side, side] may not hold either in the Critical Zone (s1, s2, f1, f2) of the Model.

Property:
The sum of the angles of a triangle can be:
- 180 degrees, if all its vertexes A, B, C are lying, for example, in delta1-f1;
- strictly less than 180 degrees [ any value in the interval (0, 180) ],
for example:
let R, T be two points in delta2-f2 such that Q does not l1e in RT, and S another point on s2;
then the triangle SRT has <(SR, ST) congruent to 0 because SR and ST have an infinite number of common points (the segment SQ), and <QTR + <TRQ congruent to 180-<TQR [ by construction we may vary <TQR in the interval (0, 180) ];
- even 0 degree!
let A be a point in delta1-f1, B a point in delta2-f2, and C a point on s3, very close to P;
then ABC is a non-degenerate triangle (because its vertexes are non-collinear), but <(A-P-s1-Q-B, A-P-s3-C)
= <(B-Q-s1-P-A, B-Q-s1-P-s3-C) = <(C-s3-P-A, C-s3-P-s1-Q-B) = 0
(one considers the length C-s3-P-s1-Q-B strictly less than C-s3-B);
the area of this triangle is also 0 !
- more than 180 degrees,
for example:
let A, B be two points in delta1-f1, such that
<PAB + <PBA + <(s1, s2; in Q) is strictly greater than 180 degrees;
then the triangle ABQ, formed by the intersection of the lines A-P-s2-Q, Q-s1-P-B, AB will have the sum of its angles strictly greater than 180 degrees.

Definition:
A circle of center M is a totality of all points A for which the segments MA are congruent to one another.

For example, if the center is Q, and the length of the segments MA is chosen greater than the length of s1, then the circle is formed by the arc of circle centered in Q, of radius MA, and lying in delta2, plus another arc of circle centered in P, of radius MA-length of s1, lying in delta1.

GROUP V. ANTI-AXIOM OF CONTINUITY (ANTI-ARCHIMEDEAN AXIOM)

Let A, B be two points. Take the points A1, A2, A3, A4, ... so that A1 lies between A and A2, A2 lies between A1 and A3, A3 lies between A2 and A4, etc. and the segments AA1, A1A2, A2A3, A3A4, ... are congruent to one another.
Then, among this series of points, not always there exists a certain point An such that B lies between A and An.

For example:
let A be a point in delta1-f1, and B a point on f1, B different from P;
on the line AB consider the points A1, A2, A3, A4, ... in between A and B, such that AA1, A1A2, A2A3, A3A4, etc. are congruent to one another;
then we find that there is no point behind B (considering the direction from A to B), because B is a limit point (the line AB ends in B).

The Bolzano's (intermediate value) theorem may not hold in the Critical Zone of the Model.

FOURTH EXAMPLE OF MODEL M7 IN RUGINA'S ORIENTATION TABLE:
**The Inconsistent System of Axioms**, and **The Contradictory Theory** (the percentage of instability is 100 **-** even... more, this is the system of chaos!)

Let (a1), (a2), ..., (an), (b) be n+1 independent axioms,
with n >= 1; and let (b') be another axiom contradictory to (b).
We construct a system of n+2 axioms:
[I ] (a1), (a2), ..., (an), (b), (b')
which is inconsistent. But this system may be shared into two consistent systems of independent axioms
[C ] (a1), (a2), ..., (an), (b),
and
[C'] (a1), (a2), ..., (an), (b').
We also consider the partial system of independent axioms
[P ] (a1), (a2), ... (an).
Developing [P ], we find many propositions (theorems, lemmas)
(p1), (p2), ..., (pm),
by combinations of its axioms.
Developing [C ], we find all propositions of [P ]
(p1), (p2), ..., (pm),
resulted by combinations of (a1), (a2), ..., (an),

plus other propositions

(r1), (r2), ..., (rt),

resulted by combinations of (b) with any of (a1), (a2), ..., (an).

Similarly for [C'], we find the propositions of [P ]

(p1), (p2), ..., (pm),

plus other propositions

(r'1), (r'2), ..., (r't),

resulted by combinations of (b') with any of (a1), (a2), ..., (an), where (r'1) is an axiom contradictory to (r1), and so on.

Now, developing [I ], we'll find all the previous resulted propositions:

(p1), (p2), ..., (pm),
(r1), (r2), ..., (rt),
(r'1), (r'2), ..., (r't).

Therefore, [I ] is equivalent to [C ] reunited to [C'].

From one pair of contradictory propositions {(b) and (b')} in its beginning, [I ] adds t more such pairs, where t >= 1, {(r1) and (r'1), ..., (rt) and (r't)}, after a complete step.

The further we go, the more pairs of contradictory propositions are accumulating in [I ].

It is interesting to study the case when n = 0.

Why do people avoid thinking about the CONTRADICTORY THEORY ?

As you know, nature is not perfect:

and opposite phenomena occur together,

and opposite ideas are simultaneously asserted and, ironically, proved that both of them are true! How is that possible? ...

A statement may be true in a referential system, but false in another one. The truth is subjective. The proof is relative.

( In philosophy there is a theory:

that "knowledge is relative to the mind, or things can be known only through their effects on the mind, and consequently there can be no knowledge of reality as it is in itself", called "the Relativity of Knowledge";

<Webster's New World Dictionary of American English>, Third College Edition, Cleveland & New York, Simon & Schuster Inc., Editors: Victoria Neufeldt, David B. Guralnik, 1988, p. 1133. )

You know? ... sometimes is good to be wrong!